\documentclass[reqno]{amsart}
\usepackage{amssymb,amsmath,amsthm,amstext,amsfonts}
\usepackage[dvips]{graphicx}
\usepackage{psfrag}
\usepackage{url}
\usepackage{amsmath,amstext,amsthm,amsfonts}
\usepackage{color}
\usepackage{xcolor} 

\usepackage[colorlinks=true, linkcolor=blue, urlcolor=red, citecolor=blue, hyperindex]{hyperref}

\makeatletter \@addtoreset{equation}{section} \makeatother

\renewcommand\thetable{\thesection.\@arabic\c@table}

\theoremstyle{plain}
\newtheorem{maintheorem}{Theorem}

\newtheorem{maincorollary}{Corollary}

\newtheorem{mainproposition}{Proposition}

\newtheorem{theorem}{Theorem }[section]

\theoremstyle{definition} \theoremstyle{remark}

\newtheorem{example}[theorem]{Example}

\newtheorem{problem}{Problem}

\newtheorem{question}{Question}
\newtheorem{conjecture}{Conjecture}

\newcommand{\Sc}{\mathbb{S}}

\newcommand{\N}{\mathbb{N}}
\newcommand{\R}{\mathbb{R}}

\newcommand{\topp}{\operatorname{top}}

\newcommand{\htop}{h_{\topp}}

\begin{document}

\title{Phase transitions for surface diffeomorphisms}

\author{ Thiago Bomfim and Paulo Varandas}

\address{Thiago Bomfim, Departamento de Matem\'atica, Universidade Federal da Bahia\\
Av. Ademar de Barros s/n, 40170-110 Salvador, Brazil.}
\email{tbnunes@ufba.br}
\urladdr{https://sites.google.com/site/homepageofthiagobomfim/}

\address{Paulo Varandas,  CMUP, Faculdade de Ci\^encias, Universidade do Porto, Rua do Campo Alegre s/n, 4169-007 Porto, Portugal \& Departamento de Matem\'atica, Universidade Federal da Bahia\\
Av. Ademar de Barros s/n, 40170-110 Salvador, Brazil}
\email{paulo.varandas@ufba.br}
\urladdr{https://sites.google.com/view/paulovarandas/}

\date{\today}

\begin{abstract}

In this paper we consider $C^1$ surface diffeomorphisms and study the existence of phase transitions, here expressed by 
the 
non-analiticity of the pressure function associated to smooth and geometric-type potentials.  
We prove that the space of $C^1$-surface diffeomorphisms admitting phase transitions is a $C^1$-Baire generic subset of the space of non-Anosov diffeomorphisms. In particular, if $S$ is a compact surface which is not homeomorphic to the 2-torus then a $C^1$-generic diffeomorphism on $S$ 
has phase transitions. We obtain similar statements in the context of $C^1$--volume preserving diffeomorphisms. Finally, we prove that a $C^2$-surface diffeomorphism exhibiting a dominated splitting admits phase transitions if and only if 
has some non-hyperbolic periodic point. 
\end{abstract}


\maketitle

\section{Introduction}

The notion of uniform hyperbolicity, coined by Smale \cite{Sm67} in the late sixties, draw the attention of many researchers and influenced the way of studying differentiable dynamical systems. In the context of non-invertible maps, hyperbolic dynamics correspond to hyperbolic endomorphism (see \cite{P76}) and uniformly expanding maps.
After the pioneering work by Sinai, Bowen and Ruelle \cite{Bow75, BR75,S72} and the advent of Markov partitions, the thermodynamic formalism of uniformly hyperbolic dynamical systems, both maps and flows, is nowadays well understood. 
Most generally, given a continuous map $f$ acting on a compact metric space $X$  and a continuous potential $\phi : X \rightarrow \R$  the variational principle for the topological pressure asserts that the topological pressure $P_{\text{top}}(f, \phi)$ of $f$ with respect to $\phi$ satisfies
\begin{equation}\label{eq:VP}
P_{\text{top}}(f, \phi) = \sup\Big\{h_{\mu}(f) + \int\phi \, d\mu : \mu \text{ is an } \text{f-invariant probability measure}\Big\},
\end{equation}
where $h_{\mu}(f)$ denotes the Kolmogorov-Sinai metric entropy of $\mu$ (see e.g. \cite{W82}).  
Whenever $\phi \equiv 0$ the topological pressure $P_{\text{top}}(f, 0)$ coincides with the topological entropy $h_{\text{top}}(f)$ of $f$, which is one of the most important topological invariants and measurements of chaotic behavior in dynamical systems. 
In case $f$ has finite topological entropy, for each continuous potential $\phi: X\to \mathbb R$ the pressure function $\R \ni t \mapsto P_{\text{top}}(f, t\phi)$ is a continuous and convex function and so,  using Aleksandrov's theorem, it is twice differentiable Lebesgue almost everywhere.
An equilibrium state $\mu_{\phi}$ for $f$ with respect to $\phi$ is an $f$-invariant probability measure which attains
the supremum in ~\eqref{eq:VP}.  
In the specific context of Maneville-Pomeau maps $f$, the non-hyperbolicity of the map is detected by a phase transition associated to the geometric potential $-\log|f'|$ (cf. \cite{Lo93,PS92}). 
 Similarly, in the context of $C^2$-partially hyperbolic diffeomorphisms admitting a $Df$-invariant splitting $TS=E^u\oplus E^{c}\oplus E^{s}$, where $E^{u}$ is the unstable subbundle and $E^{c}$ is the central subbundle, it is natural to understand the existence of equilibrium states associated to the potentials
$\phi = -\log |\det Df_{E^{u}}|$ or $\phi = -\log |\det Df_{E^{cu }}|$ (we refer the reader to \cite{LY85} for more details). 
We will say that $\phi: X\to\mathbb R$ is a \emph{geometric-type potential} for the dynamical system $f$ 
if $\phi(x)=-\log \beta(x)$ for every $x\in X$, 
where $\beta(\cdot)$ is continuous as a function of the derivative $Df(\cdot)$.
In the context of uniformly hyperbolic maps (we refer the reader to Subsection~\ref{sec:defhyp} for the definition) one has much stronger properties, among which we highlight that: 
\begin{enumerate}
\item[(P1)] every continuous potential $\phi$ admits at least one equilibrium state;
\item[(P2)] each H\"older continuous potential $\phi$ has a unique equilibrium state $\mu_\phi$;
\item[(P3)] for each H\"older continuous potentials $\phi,\psi$ the pressure function 
$$\R \ni t \mapsto P_{\text{top}}(f, \phi+t\psi)$$ is real analytic;
\item[(P4)] for every H\"older continuous potentials $\phi,\psi$ the pressure function 
$$\R \ni t \mapsto P_{\text{top}}(f, \phi+t\psi)$$ is $C^1$-differentiable.
\end{enumerate}
Let us comment on these properties. Property (P1) is a consequence of the variational principle together with the fact that the Kolmogorov-Sinai measure theoretical entropy of expansive maps is upper semicontinuous (see e.g. \cite{W82}). Moreover, property (P4) is equivalent to property (P2) by the work of Walters \cite{W92}.
Furthermore, it is worth mentioning the uniqueness of equilibrium states requires further regularity of the potentials rather than continuity (cf. \cite{Hof77}). 
\smallskip

In both the Statistical Physics and Dynamical Systems communities one often says that a phase transition occurs whenever there exists an abrupt change in the properties of the dynamical system. For that reason, the notion of phase transition  
has been used in a broad sense to describe the failure of one of the properties (P1)-(P4) listed above. 
In this broad terminology, phase transitions have been shown for several maps beyond the context of uniform hyperbolicity  (see e.g. 
\cite{CRL13, CRL15, CRL19, DGR14, IRV18, Lo93,PS92} just to mention a few). 
Following \cite{BC21}, we say that $f$ has a \emph{phase transition} with respect to the potential $\phi: X\to \mathbb R$ 
if the topological pressure function
$$
\R \ni t \mapsto P_{\text{top}}(f, t\phi)
$$
is not analytic.
In \cite[Problem~A]{BC21}, Bomfim and Carneiro asked whether it is possible to describe the mechanisms responsible for the existence of 
phase transitions for $C^{\,1}$-local diffeomorphisms with positive topological entropy and H\"older continuous potentials, and 
obtained a 
fine description of the pressure function for local diffeomorphisms on the circle. 
%
%
%
%
%
%
%
%
%
The key idea is that, in the context of circle maps the lack of expansion determines the existence of  phase transitions
and, consequently, such a $C^{1+\alpha}$-local diffeomorphism has a phase transition with respect to a H\"older continuous potential if and only if 
the same holds 
with respect to the potential $-\log |f'|$.
A similar phenomenon fails to occur in higher dimension (cf. Example~\ref{example1}). 
The methods used in \cite{BC21} were extended to describe phase transitions for 
partially hyperbolic skew-products with one dimensional center direction and the geometric-type potential which describes the central Lyapunov exponent
(cf. \cite{BCF22} for the precise statements). 
However, a characterization of the set of smooth dynamical systems 
acting on an arbitrary compact Riemannian manifold and admitting phase transitions still seems out of reach.
%
%

 \smallskip
The main goal of this paper is to study the mechanisms leading to phase transitions for surface diffeomorphisms. In our context, when we say that the dynamics has phase transititons means that the topological pressure function $t\mapsto P_{\text{top}}(f, t\phi)$ is 
non-analytic  for some 
smooth or geometric-type potential $\phi$.  
We prove that a $C^1$-generic surface diffeomorphism is either transitive Anosov or 
the pressure function associated to some smooth potentials is not smooth
(cf. proof of Theorem \ref{mainthA}), 
and obtain the non-analiticity of a certain pressure function for $C^1$-generic volume preserving surface diffeomorphisms (see Corollary \ref{mainthC}).
In the dissipative setting, the proof builds over a general principle that the presence of attractors leads to the creation
of phase transitions. More precisely, if a $C^1$-diffeomorphism admits an attractor then one can build  
a smooth potential whose pressure function admits a phase transition.
The strategy differs in the volume preserving framework, where one explores the presence of indifferent periodic points for a $C^1$-dense set of 
non-Anosov diffeomorphisms, and in this case we obtain phase transitions with respect to a geometric-type potential.
In fact, similarly to the context of Manneville-Pomeau maps, such pressure functions could be smooth but are not real analytic.

In the results discussed above we address a question about the abundance of phase transitions beyond uniform hyperbolicity. 
In case of positive entropy $C^2$-surface diffeomorphisms exhibiting a dominated splitting on the non-wandering set 
we obtain stronger results, namely the characterization of the diffeomorphisms exhibiting phase transitions.
Indeed, such diffeomorphisms are either transitive Anosov or they have phase transitions (see Theorem \ref{mainthB}).
Furthermore, as a consequence of Corollary \ref{corthm0} one concludes that 
phase transitions occur if and only if the diffeomorphism has some non-hyperbolic periodic point. 
This thermodynamic characterization of uniform hyperbolicity is in parallel to the following well known results by Ma\~n\'e saying that a $C^2$ transitive one-dimensional map without critical points 
is either hyperbolic or it has some indifferent periodic point (cf. \cite{MaCMP} for more general statements).

\smallskip
This paper is organized as follows. Section~\ref{Statement of the main results} is devoted to the statement of the main results 
on  phase transitions for surface diffeomorphisms, which involve the notions of hyperbolicity, Anosov diffeomorphisms and dominated splittings.
In Section~\ref{prelim} we recall the notions of uniform and and partial hyperbolicity, Lyapunov exponents and the Margulis-Ruelle inequality.
The proofs of the main results appear in 
Section~\ref{Proofs}. Finally, in Section~\ref{stratpro} we discuss possible extensions of these results and 
include some possible directions of research and pose a number of questions. 
%
%

\section{Statement of the main results}\label{Statement of the main results}

\subsection{Dissipative surface diffeomorphisms}\label{sec:main1}

We will denote by ${S}$ a compact connected Riemannian surface. For each integer $r\geqslant 0$ and $0\leqslant\alpha \leqslant1$ let
$\mbox{Diff}^{\,r+\alpha}(S)$ denote the space of $C^r$-diffeomorphisms on $S$ whose $r^{th}$-derivative is $\alpha$-H\"older continuous, 
endowed with the $C^{\,r+\alpha}$-topology, and let $\mathbb A^{r+\alpha}(S)$ denote the (possibly empty) open subset of $\mbox{Diff}^{\,r+\alpha}(S)$ formed by Anosov diffeomorphisms on $S$ (we refer the reader to Subsection~\ref{sec:defhyp} for the definition). 
\smallskip

Throughout the paper we will simply say that a diffeomorphism has \emph{no phase transitions} 
if the pressure function associated to potentials which are either H\"older continuous or of geometric-type is not analytic.

It is well known that Anosov diffeomorphisms have no phase transitions 
with respect to H\"older continuous potentials or the geometric potential given by the logarithm of the unstable Jacobian, which in this case
is H\"older continuous as well. 
Our first result shows that phase transitions are common in the complement of the space of Anosov diffeomorphisms. More precisely:

\begin{maintheorem}\label{mainthA}
There exists a Baire generic subset $\mathcal{R} \subset \mbox{Diff}^{\,1}(S) \setminus \mathbb A^{1}(S)$ such that, if $f \in \mathcal{R}$ then 
$f$ has a phase transition. 
\end{maintheorem}

All phase transitions appearing in the context of the previous theorem correspond to lack of differentiability of the pressure function
associated to certain smooth potentials (cf. proof of Theorem~\ref{mainthA}).
The previous theorem has an important consequence for surfaces $S$ which admit no Anosov diffeomorphisms. Franks~\cite{F70}
proved that the fundamental group $\pi_1(S)$ of a surface $S$ supporting an Anosov diffeomorphims is isomorphic to $\pi_1(\mathbb T^2)\simeq \mathbb Z$
and that every Anosov diffeomorphism on a compact surface is topologically conjugated to a hyperbolic automorphism on $\mathbb T^2$.
It is a well known fact that every codimension one Anosov diffeomorphism is transitive (cf. \cite{F70} and \cite[Theorem~1.2]{N70}), hence all Anosov diffeomorphisms on surfaces are transitive. 
Together with Theorem~\ref{mainthA}, this yields the following consequence:

\begin{maincorollary}
Let $S$ be a compact surface and assume that $S$ is not homeomorphic to the torus $\mathbb T^2$. 
There exist phase transitions 
 for all diffeomorphisms in a  $C^1$-Baire generic subset of $\mbox{Diff}^{\,1}(S)$.
\end{maincorollary}

%

\smallskip
The following result illustrates that the existence of attractors is a sufficient condition for the presence and robustness of phase transitions, even in the
context of non-transitive Axiom A diffeomorphisms.

    \begin{mainproposition}\label{mainthD}
Let $M$ be a compact and connected Riemannian manifold. For each $k \in \N$ there exists an open subset $\mathcal{O}_k \subset \mbox{Diff}^{\,1}(M)$ 
formed by Axiom A diffeomorphisms such that the following property holds: there exists a smooth potential $\phi$ so that each $f \in \mathcal{O}_k$  
has exactly $k$  phase transitions with respect to $\phi$.
\end{mainproposition}
 
In fact each open subset $\mathcal O_k$ obtained in the previous proposition is a $C^1$-open neighborhood of a non-transitive Axiom A diffeomorphism,
obtained by structural stability (cf. Subsection~\ref{pprop} for more details). Therefore, a natural question is whether we can find an open subset, as in the previous theorem, provided that no element of this open be an Axiom A diffeomorphism. 

\smallskip
The previous Theorem~\ref{mainthA} characterizes phase transitions for typical $C^1$-diffeo\-mor\-phisms, meaning a $C^1$-Baire generic subset of surface diffeomorphisms. Nonetheless, it does not characterize in such generality the geometric aspects which are the mechanisms leading to the creation of phase transitions. The next result characterizes the diffeomorphisms having phase transitions among those whose non-wandering 
$\Omega(f)$ admits a dominated splitting (cf. Subsection~\ref{sec:defhyp} for definitions).


 \begin{maintheorem}\label{mainthB}
 Let $r\geqslant 1$, $0\leqslant\alpha\leqslant1$ be such that $r+\alpha >1$.  
 Assume that $f \in \mbox{Diff}^{\,r+\alpha}(S)$ is such that $f_{|\Omega(f)}$ has a dominated splitting $T_{\Omega(f)} M = E \oplus F$ and $h_{\text{top}}(f) > 0$. 
Then, either:
\begin{itemize}
\item[(i)] $f$ is an Anosov diffeomorphism, or
\item[(ii)] $f$ has a phase transition.
\end{itemize}
 \end{maintheorem}


We observe that the positive entropy assumption in Theorem~\ref{mainthB} is necessary, and that it is automatically satisfied in case the non-wandering set $\Omega(f)$ is the entire surface $S$
(cf. Subsection~\ref{sec:defhyp}).
Moreover, using \cite{Ka80,PSa09}, one knows that if a diffeomorphism $f \in \mbox{Diff}^{2}(S)$ is such that $f_{|\Omega(f)}$ has a dominated splitting then the positive entropy assumption $h_{\text{top}}(f) > 0$ is equivalent the existence of a sequence of periodic points with unbounded periods.
Furthermore, we state a stronger result for $C^2$ surface diffeomorphisms admitting a dominated splitting, that  
every non-Anosov diffeomorphism has phase transitions. 


\begin{maincorollary}\label{corthm0}
Let $f\in \text{Diff}^{\,2}(S)$ be a transitive 
diffeomorphism on a compact surface $S$ admitting a dominated splitting $TS=E \oplus F$. The following are equivalent:
\begin{enumerate}
\item $f$ all periodic points are hyperbolic
\item $f$ has no phase transitions.
\item $f$ is an Anosov diffeomorphism
\end{enumerate}
\end{maincorollary}

Let us observe that 
if $S$ is not homeomorphic to the torus $\mathbb T^2$, the latter ensures that every volume preserving transitive diffeomorphism with a dominated splitting has an indifferent periodic point and it has phase transitions.

\medskip
Finally, the following example, inspired from \cite{BC21}, illustrates that the invertibility assumption is crucial in Theorem~\ref{mainthB}.

\begin{example}\label{example1} 
\emph{
Let $f: \mathbb S^1 \to \mathbb S^1$ be the doubling map given by $f(x)=2x \,(\text{ mod } 1)$, 
$R_{\alpha} : \Sc^{\,1}\to \Sc^{\,1}$ be an irrational rotation given by $R_{\alpha}(y)=y + \alpha \,(\text{ mod } 1)$ with $\alpha\in \mathbb R\setminus \mathbb Q$ and $F : \Sc^{\,1} \times \Sc^{\,1} \to \Sc^{\,1} \times \Sc^{\,1}$ denote the local diffeomorphism obtained as the product map
$F(x,y)=(f(x),R_\alpha(y))$ for every $(x,y)\in \Sc^{\,1} \times \Sc^{\,1}$, which has a dominated splitting inherited from the product structure and the fact that the rotation has neutral behavior.
%
Given an arbitrary H\"older continuous potential $\phi : \Sc^{\,1} \times \Sc^{\,1} \rightarrow \R$ define a H\"older continuous potential on the circle
$\tilde \phi : \Sc^{\,1} \to\mathbb R$
by
$$\tilde \phi(x)= \int \phi(x , y) \,dm(y),$$
where $m$ denotes the Lebesgue measure on the circle. As $m$ is the only $R_\alpha$-invariant probability measure then all $F$-invariant probability measures are of the form $\nu=\mu\times m$ with $\mu \in \mathcal{M}_{1}(f)$. Using the variational principle that the metric entropy of the product is additive and $h_m(R_\alpha)=0$ one concludes that
$$
P_{\text{top}}(F , \phi) = \sup_{\mu \in \mathcal{M}_{1}(f)} \Big\{h_{\mu}(f) + h_{m}(R_{\alpha}) + \int \phi(x , y)dm(y)d\mu(x) \,\Big\}= P_{\text{top}}(f , \tilde{\phi}).
$$
Since $f$ is an expanding map then $C^{\alpha}(\Sc^{\,1}  , \R) \ni \psi \mapsto P_{\text{top}}(f , \psi)$ is analytic (see e.g. \cite{BCV16,BC19})
and $\phi \mapsto \tilde \phi$ is a linear map from $C^{\alpha}(\Sc^{\,1} \times \Sc^{\,1}  , \R)$ to $C^{\alpha}(\Sc^{\,1}  , \R)$ then 
the map
$$C^{\alpha}(\Sc^{\,1} \times \Sc^{\,1}  , \R)\ni \phi \mapsto P_{\text{top}}(F , t \phi)$$ is analytic
even though $F$ is not a hyperbolic map.}
\end{example}

\subsection{Volume preserving surface diffeomorphisms}\label{sec:main2}
We also obtain analogous results in the conservative case. Given a volume measure $\omega$, we denote by $\mbox{Diff}^{\,1}_{\omega}(S)$ the space of the $C^{\,1}$-diffeomorphism that preserves $\omega$. On his space we endow the $C^{\,1}$-topology.


\begin{maintheorem}\label{mainthE}
Let $f \in \mbox{Diff}^{\,1+\alpha}_{\omega}(S)$ be such that $f$ admits a non-hyperbolic periodic point and $h_{\text{top}}(f) > 0$. Then $f$  has  phase transition.
\end{maintheorem}

Combining Theorem~\ref{mainthE} with the description of $C^1$-generic volume preserving diffeomorphisms from \cite{N78} we derive the following consequence:

\begin{maincorollary}\label{mainthC}
There is a Baire generic subset $\mathcal{R} \subset \mbox{Diff}^{\,1}_{\omega}(S)$ such that: if $f \in \mathcal{R}$ then either
$f$ is a transitive Anosov diffeomorphism or $f$ has  phase transitions.
\end{maincorollary}

\section{Preliminaries}\label{prelim}

In this section we provide some definitions and preparatory results needed for the proof of the main results.  We first recall some concepts related to uniform and partial hyperbolicity (Subsection~\ref{sec:defhyp}) and later recall an important relation between entropy and Lyapunov exponents (Subsection~\ref{RMineq}). 

\subsection{Hyperbolicity and partial hyperbolicity}\label{sec:defhyp}

Given $f\in {\rm Diff}^{\, 1}(M)$,
a $Df$-invariant splitting $TM=E\oplus F$ is
\emph{dominated} if there is an integer $k\in\mathbb{N}$ such that
\begin{equation}\label{eq:DS}
\frac{\|Df^k(x)u\|}{\|Df^k(x)w\|}<\frac{1}{2},
\end{equation}
for every $x\in M$ and every pair of unit vectors $u\in E_x$
and $w\in F_x$.
More generally, a
$Df$-invariant splitting $TM=E_1\oplus\cdots\oplus E_k$ is \emph{dominated} if
$(E_1\oplus\cdots\oplus E_l)\oplus(E_{l+1}\oplus\cdots\oplus E_k)$ is a dominated splitting, for every $1\le l\le k-1$.

\smallskip
A $Df$-invariant bundle $E$ is \emph{uniformly contracting} (resp. \emph{uniformly expanding}) if
there are constants $C>0$ and $0<\lambda<1$ such that 
$$
\|Df^n(x)v\|\le C\lambda^n\| v\| \qquad (\,\text{resp}.\; \| Df^{-n}(x)v\|\le
C\lambda^n\| v\|)
$$ 
for every $n\ge 1$, every $x\in M$ and $v\in E_x$.

\smallskip
Given a compact $f$-invariant subset $\Lambda\subset M$, we say that $\Lambda$ is \emph{partially hyperbolic} (resp. \emph{strongly partially hyperbolic})
if there is a $Df$-invariant splitting 
$$T_\Lambda M=E^s\oplus E^c\oplus E^u$$ 
such that
$E^s$ and $E^u$ are uniformly contracting and uniformly expanding respectively,
and at least one of them is (resp. both of them are) not trivial.
We say that $E^c$ is the central direction of the splitting.
The set $\Lambda$ is called \emph{hyperbolic} if it is strongly partially hyperbolic and $E^c$ is trivial.
In the special case that $\Lambda=M$ is a hyperbolic set we say that $f$ is an \emph{Anosov diffeomorphism}.

\smallskip
A diffeomorphism $f\in {\rm Diff}^{\, 1}(M)$ is called \emph{Axiom A} if $\Omega(f)$ is a hyperbolic set and the set of periodic points of $f$ is dense in the non-wandering set.  Moreover, given a compact $f$-invariant subset $\Lambda\subset M$, we say that $f\mid_\Lambda$ is \textit{transitive}
if there is $x\in \Lambda$ whose orbit is dense in $\Lambda$.

\begin{theorem}[Smale's spectral decomposition theorem (see e.g. \cite{Shub})]\label{SSDT}
If $f\in {\rm Diff}^{\, 1}(M)$ is an Axiom A diffeomorphism then there exists a unique decomposition
$$
\Omega(f) = \Lambda_1 \cup \Lambda_2 \cup \dots \cup \Lambda_k
$$
in compact, $f$-invariant, transitive and pairwise disjoint sets. 
\end{theorem}

The compact invariant sets in the statement of the previous theorem are often called \emph{hyperbolic basic pieces}. 
Recall that if $f \in \mbox{Diff}^{\,1}(M)$ and $f^{k}(p) = p$, for some $k \geqslant 1$, one says that $p$ is a \emph{hyperbolic periodic point} 
if $Df^{k}(p)$ has no eigenvalues in the unit circle.
In the modern terminology, each transitive subset $\Lambda_i$ of the non-wandering set is an homoclinic class associated to any of the hyperbolic
 periodic points contained in $\Lambda_i$(we refer the reader to \cite{Shub} for more details).
Finally, it is clear from Theorem~\ref{SSDT} that every transitive Axiom A diffeomorphism is an Anosov diffeomorphism.


\begin{theorem} (\cite{PSa09})\label{PujalsS}
Let $f\in \text{Diff}^{\,2}(S)$ and let $\Lambda\subset \Omega(f)$ be a compact invariant subset exhibiting a dominated splitting 
such that any periodic point in $\Lambda$ is hyperbolic. Then, $\Lambda=\Lambda_1\cup \Lambda_2$ where $\Lambda_1$ is a hyperbolic set and 
$\Lambda_2$ consists of a finite union of periodic simple closed curves $C_1,\dots C_n$ normally hyperbolic, and such that $f^{m_i} : C_i \to C_i$ is conjugated to an irrational rotation ($m_i$ denotes the period of the curve $C_i$).
\end{theorem}


 
\subsection{Lyapunov exponents and entropy}\label{RMineq}

Let us recall briefly the notions of Lyapunov exponents and the Margulis-Ruelle inequality relating with the measure theoretical Kolmogorov-Sinai entropy
with positive Lyapunov exponents. 

\smallskip
Given an $f$-invariant and ergodic probability measure $\mu$, by
the Oseledets' theorem~\cite{Ose68}, there exists an $f$-invariant and full $\mu$-measure set $\mathcal{R} (\mu)\subset M$, an integer $k\ge 1$
and real numbers (\emph{Lyapunov exponents})
$$
\lambda_1 \left(f, \mu\right)>\cdots >\lambda_k\left(f, \mu\right),
$$
and for every $x\in \mathcal{R} (\mu)$ there exists a splitting $T_x M =E^{1}_{x}\oplus \cdots \oplus E^{k}_{x}$ (\emph{Oseledets subspaces}) such that
$Df(x)E^{i}_{x}=E^{i}_{f(x)}$ for every $1\leqslant i \leqslant k$ and
\begin{displaymath}
\lambda _i(f, \mu) =\lim _{n\rightarrow \infty} \dfrac{1}{n}\log \| Df^n(x)v\| 
\end{displaymath}
for every $v\in E^{i}_{x}\setminus\{0\}$ and $1\leqslant i \leqslant k$. The dimension of $E^{i}_{x}$ is called the \emph{multiplicity} of $\lambda_i(\mu)$
and define the sum of positive Lyapunov exponents, considering multiplicity, defined as
$$
\lambda_+(f,\mu)=\sum_{\lambda_j(f,\mu)>0} \;  \lambda_j(f,\mu) \cdot \dim E^j_x.
$$
The sum of non--positive Lyapunov exponents $\lambda_-(f,\mu)$ is defined analogously.
If $f$ is volume preserving, i.e. $\mu$ is the volume measure, then sum of all Lyapunov exponents is zero (cf. \cite{Ose68}). 

\smallskip
We will need the following well known relation between entropy and positive Lyapunov exponents:

\begin{theorem}[Margulis-Ruelle inequality, \cite{Rue78}]\label{Marg}
Let $f : M \rightarrow M$ be a $C^{\,1}$-local diffeomorphism that preserves an $f$-invariant and ergodic probability $\mu$. Then 
$$h_\mu(f)\leqslant  \lambda_+(f,\mu)$$
\end{theorem}

\smallskip

In the general context of continuous maps defined over a compact metric space, the variational principle for the pressure asserts that
\begin{equation}\label{eqvp}
P_{\text{top}}(f,\phi) =\sup \Big\{ h_\mu(f) + \int \phi\, d\mu \colon \mu\in \mathcal M_1(f)\Big\}
\end{equation}
for every continuous potential $\phi$, 
where $P_{\text{top}}(f)$ stands for the topological entropy of $f$ (see e.g. \cite{W82}). The topological entropy $h_{\text{top}}(f)$
of $f$ coincides with $P_{\text{top}}(f,0)$
In the special case that $f$ is a $C^1$-diffeomorphism on a compact surface, if it admits a global dominated splitting $TM=E\oplus F$ then $h_{\text{top}}(f)>0$ is strictly positive, as it is bounded below by the volume growth of the stronger subbundle (cf. \cite[Proposition~2]{SX} 
for the precise statement).

\color{black}




\section{Proof of the main results}\label{Proofs}

\subsection{Proof of the Theorem \ref{mainthA}}

In \cite{M82}, Ma\~n\'e proved that there exists a $C^{\,1}$-Baire generic subset $\mathfrak R\subset \text{Diff}^{\, 1}(S)$ of diffeomorphisms 
such that any $f\in \mathfrak R$ is either an Axiom A diffeomorphism or it has infinitely many sinks or sources. 
Recall that $P_{\text{top}}(f , \phi) = P_{\text{top}}(f^{-1} , \phi)$ for every continuous potential $\phi$, hence we can, if necessary, replace the diffeomorphism $f$ by its inverse $f^{-1}$. 
Fixing  $f\in \mathfrak R$, there are two cases to consider.

\medskip
\noindent \textit{Case 1:} $f$ has two sinks (or sources)
\smallskip

\noindent  Replacing $f$ by $f^{-1}$, if necessary, assume that $p_{1}, p_{2} \in S$ are two sinks for $f$ and that $k\ge 1$ is a common period. Thus, 
$p_{1}, p_{2}$ are attracting fixed points for $f^k$ and there exist
open subsets $A_{1,i}, A_{2,i} \subset S$ such that, for each $i = 1, 2$:

\begin{itemize}
\item  $p_{i} \in {A_{1,i}}$;
\item  $\overline{A_{1,i}} \subset A_{2,i}$ are contained in the 
	basin of attraction $B(p_i)$ of $p_i$; 
    \item $f^{k}(A_{2,i}) \subset \overline{A_{1,i}}$; 
    \item $\sup_{x \in \overline{A_{2,i}}}\|Df^{k}(x)\| < 1$
\end{itemize}
and 
$$
\Big(\bigcup_{j = 0}^{k-1}f^{j}(A_{2,1}) \Big) \cap \Big(\bigcup_{j = 0}^{k-1}f^{j}(A_{2,2}) \Big) = \emptyset.
$$
Using bump functions, there exists a $C^{\infty}$ potential $\phi : S \rightarrow \R$ so that:
\begin{itemize}
\item $\phi_{|f^{j}(A_{1,1})} \equiv \frac{1}{k},$ for all $j = 0, \dots, k-1$;
\item $\phi_{|f^{j}(A_{1,2})} \equiv \frac{1}{2k},$ for all $j = 0, \dots, k-1$;
\item $\phi_{|S \setminus \bigcup_{j = 0}^{k-1}f^{j}(A_{2,1}\cup A_{2,2})} \equiv 0.$
\end{itemize}
By construction 
\begin{equation}\label{eqchoicephi}
\lim_{n\to\infty}\frac{1}{n}\sum_{i = 1}^{n-1}\phi(f^{i}(x)) = 
\begin{cases}
\begin{array}{ll}
1 & , \text{if \,} \, x \in \bigcup_{j=0}^{k-1} f^j(B(p_1)) \\
\frac12 & , \text{if \,} \, x \in \bigcup_{j=0}^{k-1} f^j(B(p_2)) \\
0 & , \text{otherwise.} 
\end{array}
\end{cases}
\end{equation}

Now, observe that Birkhoff's Ergodic Theorem ensures that for each $\mu \in \mathcal{M}_{\text{erg}}(f)$ there exists $x \in S$ such that 
$\lim_{n\to\infty} \frac{1}{n}\sum_{i=1}^{n-1}\delta_{f^{i}(x)} = \mu$. In particular, using the ~\eqref{eqchoicephi} any $f$-invariant ergodic probability measure $\mu$ which is not supported on the periodic attractors satisfies 
$\int \phi \, d\mu = 0$ and 
\begin{equation}\label{eqchoicephiP}
P_{\mu}(f , t\phi) := h_{\mu}(f) + t\int \phi \,d\mu = h_{\mu}(f).
\end{equation}
As $P_{\delta_{p_{i}}}(f  , t\phi) = \frac{1}{i}t$ for every $i=1,2$, using the variational principle ~\eqref{eqvp} it follows that 
$$
P_{\text{top}}(f , t\phi) =  \max\Big\{h_{\text{top}}(f) , t \Big\} \qquad \text{for every}\; t\in \mathbb R.
$$
Thus $t \mapsto P_{\text{top}}(f , t\phi)$ is not differentiable at $t = h_{\text{top}}(f)$. 

%

\bigskip
\noindent \textit{Case 2:} $f$ is a non-transitive Axiom A diffeomorphism
\smallskip

\noindent 
By the spectral decomposition theorem (Theorem~\ref{SSDT}) one can write $\Omega(f)=\bigcup_{i=1}^k \Lambda_i$ where 
$\Lambda_{1}, \dots, \Lambda_{k}$ are pairwise disjoint hyperbolic basic sets. By the Poincar\'e recurrence theorem, 
$\mu(\Omega(f))=1$ for every $\mu \in \mathcal{M}_{\text{erg}}(f)$ and, consequently, $h_{\text{top}}(f) = h_{\text{top}}(f\mid_{\Omega(f)})
= \max \{  h_{\text{top}}(f\mid_{\Lambda_i}) \colon 1\leqslant i \leqslant k\}$.
Moreover, there exists a $C^\infty$-potential $\phi$ such that 
$\phi\mid_{\Lambda_1} \equiv 1$ and $\phi\mid_{\Lambda_j} \equiv -1$ for all $j = 2, \dots, k$ and, consequently, 
$$
P_{\text{top}}(f , t\phi) =  \max\Big\{ h_{\text{top}}(f\mid_{\Lambda_{1}}) + t \; , \; \max_{2\leqslant i \leqslant k} \big\{ h_{\text{top}}(f\mid_{\Lambda_{i}}) - t  \big\} \Big\}
$$
for every $t\in \mathbb R$. This pressure function is not differentiable at 
$$t=\frac12(\max_{2\leqslant i \leqslant k} \big\{ h_{\text{top}}(f\mid_{\Lambda_{i}})
-h_{\text{top}}(f\mid_{\Lambda_{1}})).$$
This concludes the proof of the Theorems \ref{mainthA}. \hfill $\square$ 

\subsection{Proof of Proposition~\ref{mainthD}}

Let $f \in \mbox{Diff}^{\,1}(S)$ be such that there exists $\Lambda_{1}, \dots, \Lambda_{k}$ disjoints hyperbolic attractors for $f$ or $f^{-1}$, with $k\geqslant 1$. In case $k >1$ we can perform a surgery and modify $f$ in such a way that 
$$h_{\text{top}}(f_{|\Lambda_{1}}) > h_{\text{top}}(f_{|\Lambda_{2}}) > \dots > h_{\text{top}}(f_{|\Lambda_{k}})$$
and
$$h_{\text{top}}(f) - h_{\text{top}}(f_{|\Lambda_{1}}) < h_{\text{top}}(f_{|\Lambda_{1}}) - h_{\text{top}}(f_{|\Lambda_{2}}) < \dots < h_{\text{top}}(f_{|\Lambda_{k-1}}) - h_{\text{top}}(f_{|\Lambda_{k}}).$$ By structural stability (see e.g. \cite{Shub}),
there exists $\mathcal{A} \subset \mbox{Diff}^{\,1}(S)$ an open neighborhood of $f$ and $\Lambda_{i} \subset  \overline{A_{1,i}} \subset A_{2,i}$ with $A_{1,i}, A_{2,i}$ open subsets of $S$, for $i = 1, \dots, k$, such that every $g \in \mathcal{A}$ satisfies:

If $\Lambda_{i}$ is attractor for $f$ then
\begin{itemize}
    \item $g(A_{2,i}) \subset \overline{A_{1,i}}$;
    \item $\Lambda_{i , g} := \bigcap_{n = 0}^{+\infty}g^{n}(A_{2,i})$ is a hyperbolic attractor for $g$;
    \item $h_{\text{top}}(g_{|\Lambda_{i,g}}) = h_{\text{top}}(f_{|\Lambda_{i}})$.
\end{itemize}

If $\Lambda_{i}$ is attractor for $f^{-1}$ then
\begin{itemize}
    \item $g^{-1}(A_{2,i}) \subset \overline{A_{1,i}}$;
    \item $\Lambda_{i , g} := \bigcap_{n = 0}^{+\infty}g^{-n}(A_{2,i})$ is a hyperbolic attractor for $g^{-1}$;
    \item $h_{\text{top}}(g_{|\Lambda_{i,g}}) = h_{\text{top}}(f_{|\Lambda_{i}})$.
\end{itemize}

Define $\phi : S \rightarrow \R$ a $C^{\infty}$ potential such that:

\begin{itemize}
\item $\phi_{|A_{1,i}} \equiv i$ for all $i=1,\dots, k$;
\item $\phi_{|S \setminus \bigcup_{i = 1  }^{k+1}A_{2,i}} \equiv 0.$
\end{itemize}

Note that $\lim_{n\to\infty}\frac{1}{n}\sum_{i = 1}^{n-1}\phi(f^{i}(x)) = 0, 1, \dots, $ or $k$, for all $x \in S$. 
Given $\mu \in \mathcal{M}_{\text{erg}}(f)$, by Birkhoff's Ergodic Theorem, there exists $x \in S$ such that 
$$\lim_{n\to\infty} \frac{1}{n}\sum_{i=1}^{n-1}\delta_{f^{i}(x)} = \mu.$$ In particular, $\mu \in \mathcal{M}_{\text{erg}}(g_{|\Lambda_{i , g}})$ for $i = 1, \dots k$ or $\int \phi \,d\mu = 0$.
If $\int \phi \,d\mu = 0$ then $P_{\mu}(g , t\phi) = h_{\mu}(f).$
On other hand, if  $\mu \in \mathcal{M}_{\text{erg}}(g_{|\Lambda_{i , g}})$ then $P_{\mu}(g , t\phi) = h_{\mu} + ti$.
Therefore,
$$
P_{\text{top}}(g , t\phi) =  \max\{h^{\ast} , h_{\text{top}}(f_{|\Lambda_{1}}) + t, \dots,  h_{\text{top}}(f_{|\Lambda_{k+1}}) + t(k+1)\},
$$ where $h^{\ast} := \sup \{h_{\mu}(f) : \mu \in \mathcal{M}_{\text{erg}}(f) \text{ and } \int \phi \,d\mu = 0\}$. Since that $f$ is a non-transitive Axiom A diffeomorphism then $h^{\ast}$ is well defined and it is strictly positive.

Observe that $t \mapsto P_{\text{top}}(g , t\phi)$ is not differentiable at $t = h^{\ast} - h_{\text{top}}(f_{|\Lambda_{1}})$. In case $k > 1$ then the pressure function $t \mapsto P_{\text{top}}(g , t\phi)$ will not be differentiable at each parameter $t_{i} = h_{\text{top}}(f_{|\Lambda_{i}}) -  h_{\text{top}}(f_{|\Lambda_{i+1}})$ as well,  where $i = 1, \dots, k-1$.

\subsection{Proof of Theorem \ref{mainthB}}\label{pprop}

 Assume that $f \in \mbox{Diff}^{\,r+\alpha}(S)$  ($r+\alpha >1$) is such that $f_{|\Omega(f)}$ has a dominated splitting $T_{\Omega(f)} M = E \oplus F$,
 that $h_{\text{top}}(f) =h_{\text{top}}(f\mid_{\Omega(f)})  > 0$ and that $f$ is not an Anosov diffeomorphism. 





\medskip
If $f$ 
is a non-transitive Axiom A diffeomorphism then $f$ has  phase transition with respect to some smooth potential
(recall {Case 2} in the proof of Theorem~\ref{mainthA}). 
Otherwise, $f$  
has a dominated splitting $T_{\Omega(f)}S = E \oplus F$ but it is not an Axiom A diffeomorphism. 

\smallskip
We claim that either the subbundle $F$ is not uniformly contracting or the subbundle $E$ is not uniformly expanding. 
Indeed, if $F$ is uniformly contracting (resp. $E$ is uniformly expanding), using that the splitting is dominated (recall ~\eqref{eq:DS}), 
then the subbundle $E$ would be 
uniformly contracting (resp. the subbundle $F$ would be 
uniformly expanding), which would imply that $f^{-1}$ will be an expanding dynamics, in particular $f^{-1}$ is not invertible.

\smallskip
Hence we may assume, without loss of generality, that the subbundle $F$ is not uniformly contracting.
By \cite[Lemma I.5]{M88}, as $F$ is not uniformly contracting then there exists an $f$-invariant and ergodic probability 
measure $\mu \in \mathcal{M}_{\text{erg}}(f)$ so that  $\int \log \|Df_{|F}\| d\mu \geqslant 0$.
We observe that there exists an $f$-invariant and ergodic probability  $\nu \in \mathcal{M}_{\text{erg}}(f)$ with $\int \log \|Df\mid_F\| d\nu > 0$. Indeed, if this was not the case then $\int \log \|Df_{|F}\| d\eta \leqslant 0$ for every $\eta \in \mathcal{M}_{\text{erg}}(f)$ and, by the Margulis-Ruelle inequality (Theorem~\ref{Marg}), one would get
$$
0\leqslant h_\eta(f) \leqslant \int \log \|Df_{|F}\| d\eta \leqslant 0, \qquad \forall \eta \in \mathcal{M}_{\text{erg}}(f),
$$
and, consequently, $\htop(f)=\sup_{\eta \in \mathcal{M}_{\text{erg}}(f)}\,\{ h_\eta(f)\}=0$, which contradicts the positive entropy assumption.
A similar argument applies to the subbundle $E$. Define 
$$
{\lambda_{\min}}^{\!\!\!\!\!F} := \min_{\eta\in \mathcal{M}_{\text{erg}}(f)} \int \log \|Df_{|F}\| d\eta \; ,\; 
 {\lambda_{\max}}^{\!\!\!\!\!F} := \max_{\eta\in \mathcal{M}_{\text{erg}}(f)}\int \log \|Df_{|F}\| d\eta , 
$$
and analogously for the subbundle $E$.
Therefore, replacing $f$ by $f^{-1}$ if necessary, we can assume without loss of generality that either
$$
{\lambda_{\max}}^{\!\!\!\!\!F} > 0 >  {\lambda_{\min}}^{\!\!\!\!\!F}
\quad\text{ or }\quad  
{\lambda_{\max}}^{\!\!\!\!\!F} > 0 =  {\lambda_{\min}}^{\!\!\!\!\!F} \geqslant {\lambda_{\max}}^{\!\!\!\!\!E}.
$$

Now, define the potential  $\phi: S \to \mathbb R$ by $\phi(x) := -\log \|Df_{|F_{x}}(x)\|$ and the Lyapunov exponent 
${\lambda_F}(\mu) := \int \log \|Df_{|F}\| d\mu $.
Consider the non-empty convex sets 
$$
\mathcal{M}^{+}(f) := \{\mu \in\mathcal{M}_{\text{erg}}(f) : {\lambda_F}(\mu) > 0\}
$$
$$
\mathcal{M}^{-}(f) := \{\mu \in\mathcal{M}_{\text{erg}}(f) : {\lambda_F}(\mu) \leqslant 0\}
$$
and the corresponding pressure functions $P^\pm : \mathbb R \to \mathbb R$ given by
$$
P^{+}(t) := 
\sup_{\mu \in \mathcal{M}^{+}(f)} \; \Big\{h_{\mu}(f) + t\int \phi \,d\mu\Big\} 
= 
\sup_{\mu \in \mathcal{M}^{+}(f)} \; \Big\{h_{\mu}(f) + t \, \lambda_F(\mu) \Big\} 
$$
and
$$
P^{-}(t) := 
\sup_{\mu \in \mathcal{M}^{-}(f)} \Big\{h_{\mu}(f) + t\int \phi \, d\mu\Big\}
= 
\sup_{\mu \in \mathcal{M}^{-}(f)} \; \Big\{h_{\mu}(f) + t \, \lambda_F(\mu) \Big\} .
$$
Note that 
$P^{+}$ is a convex decreasing map whereas $P^{-}$ is a convex increasing map of the real line and, using the variational principle,
\begin{equation}\label{eq:aux1}
P_{\text{top}}(f , t\phi) = \max\{P^{+}(t) , P^{-}(t)\}
\end{equation}
and
\begin{equation}\label{eq:aux2}
0 \leqslant P^{+}(0) \leqslant h_{\text{top}}(f).
\end{equation}
We will use the following auxiliary estimates on the latter pressure functions:

\medskip
\noindent{\bf Claim:} \emph{$P^{-}(t) = -t{\lambda_{\min}}^{\!\!\!\!\!F}$ for all $t \geqslant 0$.}
\smallskip

\begin{proof}[Proof of the Claim]
The dominated splitting assumption implies that if $\eta \in \mathcal{M}^{-}(f)$ then $\int \log \|Df\mid_E\| d\eta \leq 0$.
Applying the Margulis-Ruelle's inequality one concludes that $0\leqslant h_{\mu}(f) \leqslant {\lambda_F}(\mu) \leqslant 0$ for all $\eta \in \mathcal{M}^{-}(f)$. This proves that $P^{-}(t) = - t {\lambda_{\min}}^{\!\!\!\!\!F}$ for all $t \geqslant 0$, as desired.
\end{proof}

\smallskip

Suppose that ${\lambda_{\min}}^{\!\!\!\!\!F} < 0$.
Note that $P(t) = P^{+}(t)$ for all $t \leqslant 0$. Moreover, by the Claim one has that
 $P^{+}(t) \leqslant  h_{\text{top}}(f) < P^{-}(t)$ for all $t > -\frac{h_{\text{top}}(f)}{{\lambda_{\min}}^{\!\!\!\!\!F}}$. Taking 
$$
t_{0} := \inf\Big\{t \in \Big[0 , -\frac{h_{\text{top}}(f)}{{\lambda_{\min}}^{\!\!\!\!\!F}}\Big] : P^{+}(t) \leqslant P^{-}(t)\Big\}
$$ 
one has that 
$$
\frac{dP}{dt}(t_0^-) \leqslant 0 \quad \text{and} \quad \frac{dP}{dt}(t_0^+)= - \lambda_{\min}^{\!F} > 0,
$$
which proves that $P(\cdot)$ is not differentiable at $t_{0}$. 

\smallskip
Suppose now that ${\lambda_{\min}}^{\!\!\!\!\!F} = 0 \geqslant {\lambda_{\max}}^{\!\!\!\!\!E}$. On the one hand, given $\mu \in \mathcal{M}_{\text{erg}}(f)$ we have $h_{\mu}(f) \leqslant \lambda_{F}(\mu)$, by the Margulis-Ruelle inequality (Theorem~\ref{Marg}). In particular, $P_{top}(f , t\phi) \leqslant 0$ for all $t\geqslant 1$.
On the other hand, for any $\eta \in \mathcal{M}_{\text{erg}}(f)$ such that $\lambda_{F}(\eta) = 0$ one obtains that $h_{\eta} =  \lambda_{F}(\eta) = 0$. This proves that $P_{top}(f , t\phi) = 0$ for all $t \geqslant 1$.
Taking 
$$
t_{0} := \inf\Big\{t \in [0 , 1] : P_{top}(f , t \phi) = 0\Big\}
$$ 
it follows that
$
P_{top}(f , t \phi) = 0$ for all $t \geqslant t_0$, and $P_{top}(f , 0\cdot\phi) = h_{top}(f) > 0$.
This proves that $t \mapsto P_{top}(f , t\phi)$ is not analytic at $t_{0}$, and completes the proof of the theorem.
\hfill $\square$

\subsection{Proof of Corollary~\ref{corthm0}}
Let $f\in \text{Diff}^{\,2}(S)$ be a transitive
diffeomorphism on a compact surface $S$ admitting a dominated splitting $TS=E \oplus F$. In particular $h_{\text{top}}(f)>0$ (cf. \cite{SX}).

Applying Theorem \ref{mainthB}, $f$ is not an Anosov diffeomorphism if and only if $f$ has a phase transition.

This proves that item (2) is equivalent to item (3). 

\smallskip
It is clear that item (3) implies item (1). Hence, in order to finish the proof of the corollary it is enough to prove 
that if all periodic points of $f$ are hyperbolic 
then $f$ is an Anosov diffeomorphism. 
Indeed, if this is the case, Theorem~\ref{PujalsS} guarantees
that $\Omega(f)=\Lambda_1\cup \Lambda_2$ where $\Lambda_1$ is a hyperbolic set and 
$\Lambda_2$ consists of a finite union of periodic simple closed curves $C_1,\dots C_n$ normally hyperbolic, and such that $f^{m_i} : C_i \to C_i$ is conjugated to an irrational rotation ($m_i$ denotes the period of the curve $C_i$). As $f$ is transitive then $\Lambda_2=\emptyset$, thus
 $\Omega(f)$ is a hyperbolic set. Thus every ergodic measure of $f$ has non-zero Lyapunov exponents. Applying the \cite[Lemma I.5]{M88}, $E$ will be uniformly contracting and $F$ will be uniformly expanding. We conclude that $f$ is Anosov.
\hfill $\square$

\subsection{Proof of the Theorem \ref{mainthE}}
Let $p \in S$ be a non-hyperbolic periodic point and $\{\lambda_1, \lambda_2\}\subset 
\mathbb C$ denote the generalized eigenvalues for $Df^k(p)$. As $p$ is not hyperbolic and $f$ is volume preserving then $|\lambda_i|=1$ for every $i=1,2$. 
 Consider the continuous
potential 
$$
\phi(x)=-\frac{1}{k}\log |\lambda(f^k , x)| , \qquad \text{ for every } \, x\in S,
$$ 
where $\lambda(f^k , x)$ denotes the eigenvalue of $Df^{k}(x)$ with largest absolute value.

We claim that $h_{\mu}(f) \leqslant - \int \phi \,d\mu,$ for every $\mu \in \mathcal{M}_{\text{erg}}(f)$. Indeed,  using Birkhoff's ergodic theorem, Kingman's sub-additive ergodic theorem and the Margulis-Ruelle inequality (recall Theorem~\ref{Marg}), one concludes that 
\begin{align*}
h_{\mu}(f)  \leqslant \lambda^{+}(f , \mu)
		 & = \lim_{n\to \infty} \frac{1}{n}\log \|Df^{n}(x)\| = \lim_{n \to +\infty}\frac{1}{nk}\log\|Df^{nk}(x)\| \\
		 & \leqslant \lim_{n \to +\infty}\frac{1}{nk}\sum_{i=0}^{n-1}\log\|Df^{k}(f^{ik}(x))\| = - \int \phi \, d\mu
\end{align*}
for each $\mu \in \mathcal{M}_{\text{erg}}(f)$ and $\mu$-almost every $x$. By the variational principle it follows that
$P_{\text{top}}(f , \phi) \leqslant 0$.
Conversely, taking the empirical measure $\mu := \frac{1}{k}\sum_{i= 0}^{k-1}\delta_{f^{i}(p)}$ one has that $h_{\mu}(f) = 0 = \int \phi \,d\mu$ and, thus, $P_{\text{top}}(f , \phi) = 0$ and $\mu$ is an equilibrium state with respect to the potential $\phi$. 

\smallskip
The topological pressure function $P(t) := P_{\text{top}}(f , t\phi)$ is a convex function of the real line such that $P(0) = h_{\text{top}}(f) > 0$. Since $f$ is 
volume preserving then $|\lambda(f , x)|\geqslant 1$ for every $x\in X$ and so, $\phi(x) \leqslant 0$ for every $x\in S$. This guarantees that
$\int \phi \, d\mu \leqslant 0$ for all $\mu \in \mathcal{M}_{\text{erg}}(f)$ and that $P(\cdot)$ is decreasing. Furthermore, $P(t) = 0$ for all $t \geqslant 1$. Defining $t_{0} : =\min\{t \in (0 , 1] : P(t) = 0\}$ we conclude that $P(\cdot)$ is not analytic at $t_{0}$. This finishes the proof of Theorem~\ref{mainthE}.
\hfill $\square$

\subsection{Proof of the Corollary \ref{mainthC}}

In \cite{N77}, Newhouse proved that there exists a Baire generic subset $\mathcal{R}_{1} \subset \mbox{Diff}^{\,1}_{\omega}(S)$ such that if $f \in \mathcal{R}_{1}$ then either $f$ is a transitive Anosov diffeomorphism or both the elliptic periodic orbits and the hyperbolic periodic orbits of $f$ 
are dense in $S$. 
Moreover, applying \cite{N78}, there exists a Baire generic subset $\mathcal{R}_{2} \subset \mbox{Diff}^{\,1}_{\omega}(S)$ such that if $f \in \mathcal{R}_{2}$ is not an Anosov diffeomorphism then, given a hyperbolic periodic point $p \in S$ with  $f^{k}(p) = p$ one has $h_{\text{top}}(f) \geqslant \frac{1}{k}\log |\lambda(f , p)|$, where $\lambda(f , p)$ is the eigenvalue of $Df^{k}(p)$ with largest absolute value. 
Altogether, if $f \in \mathcal{R}_{1}\cap \mathcal{R}_{2}$ is not an Anosov diffeomorphims  then $f$ has elliptic orbits and $h_{\text{top}}(f) > 0$. The corollary now follows as a direct consequence of Theorem \ref{mainthE}.
\hfill $\square$

\section{Further comments and questions}\label{stratpro}

\subsection{Geometric potentials}

There are certain potentials that give rise to very important invariant probabilities, such as: a.c.i.p., physical measures or SRB measures. 
In the specific context of Maneville-Pomeau maps $f$, the non-hyperbolicity of the map is detected by a phase transition associated to the geometric potential $-\log|f'|$ (cf. \cite{Lo93,PS92}). Nevertheless, in some situations the geometric-type potentials need not be regular, say H\"older continuous.  Therefore it is natural to pose the following question:

\begin{question}\label{queA}
Let $f\in \mbox{Diff}^{\,1}(S)$ be a diffeomorphism acting on a compact surface $S$. If there exists a phase transition with respect to 
H\"older continuous potentials does there exist a phase transition with respect to a geometric-type potential?
\end{question}

\subsection{Conformal diffeomorphisms}

As the reader may have realized from the arguments used in the proofs, the low dimensional assumption in Theorems~\ref{mainthA}
- ~\ref{mainthE}  is used to guarantee that the top Lyapunov exponent can be computed as a limit of a Birkhoff sum (recall e.g. the proof of Theorem~\ref{mainthE}). In particular, the results do extend in a straightforward manner to the context of conformal diffeomorphisms on a 
compact Riemannian manifold.

\subsection{Three-dimensional flows}

Non-singular vector fields on three-dimensional compact manifolds can sometimes be understood using two-dimensional maps,
as in the special case that the flows admit a global cross-section and smooth first return Poincar\'e map.
In general, the uniform hyperbolicity can be detected using the behavior of the flow in the orthogonal direction to the vector field. Indeed,
a three dimensional flow is hyperbolic if and only if the linear Poincar\'e flow is hyperbolic (cf. \cite{Doering}). Moreover, a
$C^1$-vector field either generates a hyperbolic flow or it can be $C^1$-approximated by a vector field which exhibits a homoclinic tangency or a singular cycle (cf. \cite{AR} for the definitions and statements), and the mechanisms can be used to produce periodic sinks or sources.

\begin{question}\label{queA2}
Let $X\in \mathfrak{X}^{\,1}(M)$ be a vector field on a three-dimensional compact Riemmanian manifold $M$. 
If the flow $(X_t)_{t\in\mathbb R}$ generated by $X$ is not a transitive Anosov flow does there exist arbitrary $C^1$-close vector fields which 
admit phase transitions with respect to some H\"older continuous potential or a geometric-type potential?
\end{question}

If the vector field admits a global cross section then one can answer positively Question~\ref{queA2} 
using Theorem~\ref{mainthA}. More generally, the answer to it may not necessarily make use of reduction to
piecewise smooth Poincar\'e maps acting on two-dimensional surfaces (possibly with boundary), as geometric-type
potentials like $\phi(x)=\log \|DX_t^{\wedge 2 }(x)\|$,  
where $DX_t^{\wedge 2 }(x)$ denotes the exterior power of the derivative $DX_t(x)$,
can be used to detect the uniform hyperbolicity of the vector field.

%
%

\subsection{Orders of phase transitions}

In order to obtain a dichotomy between uniform hyperbolicity and the (generic) existence of phase transitions, the latter concept
was taken in a broad sense, corresponding to either lack of differentiability of lack of analiticity of some pressure function associated
to relevant classes of potentials.   
This broad definition was used just once in the proof of Theorem~\ref{mainthB}, in the context of positive entropy diffeomorphisms 
admitting a dominated splitting, to ensure that even if the pressure function associated to geometric-type potentials is $C^1$-smooth then 
it fails to be analytic. Similarly to the context of Manneville-Pomeau maps (see e.g. \cite{Lo93}), we expect this to occur for certain diffeomorphisms on surfaces.

\begin{problem}\label{pqueA2}
Construct a smooth parameterized family $(f_\beta)_{t\in \,[0,1]}$ of $C^{1+\alpha}$-diffeomor\-phisms on a compact surface so that 
the pressure function of $f_\beta$ associated to a geometric-type potential is not smooth for every $ \beta\in [0,1)$ and it is $C^2$-smooth  
but not analytic for $\beta=1$.
\end{problem}

A natural class of diffeomorphisms to consider is the class of almost Anosov diffeomorphisms considered in \cite{HZ}, which admit
a neutral fixed point.

%
%
%
%

\subsection{Infinite phase transitions}
 Recall that there exist open sets of diffeomorphisms with arbitrarily large number of phase transitions (cf. Proposition~\ref{mainthD}).
 This suggests to ask whether one can construct a surface diffeomorphism with infinitely many phase transitions. Such constructions
 appeared e.g. for shift dynamics and continuous potentials (see \cite{KQW21} for details), but the potential is not H\"older continuous.  
 We expect that  infinite phase transitions are not possible in this differentiable context.
 
 \begin{conjecture}
 There is no diffeomorphism $f \in \mbox{Diff}^{\,1}(S)$ such that $f$ has infinite  phase transition with respect to a H\"older continuous potential.
 \end{conjecture}

%
%


\vspace{.3cm}
\subsection*{Acknowledgements}
TB was partially supported by CNPq-Brazil and Universal-18/2021. 
PV was partially supported by CMUP, which is financed by national funds through 
FCT - Funda\c c\~ao para a Ci\^encia e a Tecnologia, I.P., under the project with reference UIDB/00144/2020. PV also acknowledge financial support from the projects PTDC/MAT-PUR/29126/2017 and PTDC/MAT-PUR/4048/2021, and benefited from the grant CEECIND/03721/2017 of the Stimulus of Scientific Employment, Individual Support 2017 Call, awarded by FCT.


\bibliographystyle{alpha}

\end{document}